\documentclass[11pt]{amsart}

\usepackage{amsmath}
\usepackage{amssymb}
\usepackage{graphicx}

\theoremstyle{definition}

\numberwithin{equation}{section}

\title[Solvability of quadratic integral systems]
{Solvability of a class of systems of quadratic integral equations}

\author[Y. Chen]{Yuming Chen}
\address[Y. Chen]{Department of Mathematics, Wilfrid Laurier University,
Waterloo, Ontario, N2L 3C5, Canada}
\email{\tt ychen@wlu.ca}

\author[V. Vougalter]{Vitali Vougalter}

\address[V. Vougalter]{Department of Mathematics,  University of Toronto,
Toronto, Ontario, M5S 2E4, Canada}
\email{{\tt vitali@math.toronto.edu}}



\begin{document}

\begin{abstract}
The article is devoted to the existence of solutions of a certain system
of quadratic integral equations in $H^{1}({\mathbb R}, {\mathbb R}^{N})$. 
We show the existence of a perturbed solution by using
a fixed point technique in the Sobolev space on the real line.
\end{abstract}

\maketitle


\setcounter{equation}{0}

\section{Introduction}

\noindent
The present work deals with the existence of 
solutions of the following system of integral equations,
\begin{equation}
\label{p}
u_{m}(x)=u_{0, m}(x)+[T_{m}u_{m}(x)]
\int_{-\infty}^{\infty}K_{m}(x-y)g_{m}(u(y))dy, \quad x\in {\mathbb R},
\end{equation}
with $1\leq m\leq N$. The vector function involved in problem (\ref{p}) is
given by
\begin{equation}
\label{vf}
u(x)=(u_{1}(x), u_{2}(x), ..., u_{N}(x))^{T}.  
\end{equation}
Similarly,
\begin{equation}
\label{vf0}
u_{0}(x)=(u_{0, 1}(x), u_{0, 2}(x), ..., u_{0, N}(x))^{T}.  
\end{equation}
The exact conditions on $u_{0}(x), \ g_{m}(u)$, the linear operators
$T_{m}$ and the kernels $K_{m}(x)$ will be stated below. For each equation
of system (\ref{p}) the second term in its
right side is a product of $T_{m}u_{m}(x)$ and the integral operator
acting on the function $g_{m}(u)$, for which the sublinear growth will be
demonstrated in the proof of Theorem 1.3. further down. Thus, the system of
integral equations of this kind is called {\it quadratic}.

\noindent
The theory of 
integral equations has many useful applications in describing numerous
events and phenomena of the real world. It is caused by the fact that such
theory
is frequently applicable in various branches of mathematics, in mathematical
physics, economics, biology as well as in dealing with some real world problems.
Quadratic integral equations arise in the theories of the radiative
transfer and neutron transport, in the kinetic theory of gases, in the design
of the bandlimited signals for the binary communication using the simple
memoryless correlation detection when signals are disturbed by 
additive white Gaussian noises (see e.g. ~\cite{AC08}, ~\cite{A92}, ~\cite{E05}
and references therein).
Article ~\cite{AC08} is devoted to the
solvability of a nonlinear quadratic integral equation in the Banach space
of real functions being defined and continuous on a bounded and closed
interval using the fixed point result. Works ~\cite{AC11} and ~\cite{AC14}
deal with the existence of solutions for 
quadratic integral equations on unbounded intervals. The existence of
solutions for quadratic integral inclusions was discussed in ~\cite{AC13}.
In ~\cite{DH12} the authors study the nondecreasing solutions of a
quadratic integral equation of Urysohn-Stieltjes type.
The solvability of
quadratic integral equations in Orlicz spaces was considered in
~\cite{CM12}, ~\cite{CM15}, ~\cite{CM16}. 
The integro-differential equations, which may involve either Fredholm or
non-Fredholm operators, arise in mathematical biology when
studying systems with nonlocal consumption of resources and 
intra-specific competition (see ~\cite{EV22}, ~\cite{VV21},
~\cite{VV22} and references therein). 
The contraction argument was used
in ~\cite{V10} to evaluate the perturbation to the standing solitary wave
of the Nonlinear Schr\"odinger (NLS) equation when either the external
potential or the nonlinear term was perturbed. The analogous ideas were used
to demonstrate the persistence of pulses for some reaction-diffusion type
equations (see ~\cite{CV21}).

\noindent
Suppose that the assumption below is satisfied.

\medskip

\noindent
{\bf Assumption 1.1.}  {\it For $1\leq m\leq N$, the kernels
$K_{m}(x): {\mathbb R}\to {\mathbb R}$ do not vanish identically on the real
line, so that $K_{m}(x)\in W^{1, 1}({\mathbb R})$.
The functions
$u_{0, m}(x): {\mathbb R}\to {\mathbb R}, \ u_{0, m}(x)\in H^{1}({\mathbb R})$
and $u_{0, m}(x)$ are nontrivial for some $1\leq m\leq N$.
We also suppose that the linear operators
$T_{m}: H^{1}({\mathbb R})\to H^{1}({\mathbb R})$ are bounded, so that their
norms $0<\|T_{m}\|<\infty$.}  

\bigskip

\noindent
Here we use the Sobolev space
\begin{equation}
\label{H1}  
H^{1}({\mathbb R}):=\big\{\phi(x):{\mathbb R}\to {\mathbb {\mathbb R}} \ | \
\phi(x)\in L^{2}({\mathbb R}), \ \frac{d\phi}{dx}\in
L^{2}({\mathbb R}) \big \}.
\end{equation}
It is equipped with the norm
\begin{equation}
\label{n}
\|\phi\|_{H^{1}({\mathbb R})}^{2}:=\|\phi\|_{L^{2}({\mathbb R})}^{2}+
\Big\|\frac{d\phi}{dx}\Big\|_{L^{2}({\mathbb R})}^{2}.
\end{equation}
For a vector function (\ref{vf}) we use the norm
\begin{equation}
\label{nvf}
\|u\|_{H^{1}({\mathbb R}, {\mathbb R}^{N})}^{2}:=\|u\|_{L^{2}({\mathbb R}, {\mathbb R}^{N})}^{2}+
\sum_{m=1}^{N}\Big\|\frac{du_{m}}{dx}\Big\|_{L^{2}({\mathbb R})}^{2},
\end{equation}
with
\begin{equation}
\label{nvf2}
\|u\|_{L^{2}({\mathbb R}, {\mathbb R}^{N})}^{2}:=\sum_{m=1}^{N}\|u_{m}\|_{L^{2}({\mathbb R})}^{2}.
\end{equation}
Another norm relevant to our argument is given by
\begin{equation}
\label{11}
\|K_{m}\|_{W^{1, 1}({\mathbb R})}:=\|K_{m}\|_{L^{1}({\mathbb R})}+
\Big\|\frac{dK_{m}}{dx}\Big\|_{L^{1}({\mathbb R})}.
\end{equation}
We define
\begin{equation}
\label{Q}
Q:=\sqrt{\sum_{m=1}^{N}\|T_{m}\|^{2}\|K_{m}\|_{W^{1, 1}({\mathbb R})}^{2}},  
\end{equation}
such that $0<Q<\infty$.

\noindent
Let the functions $V_{m}(x): {\mathbb R}\to {\mathbb R}$ be nontrivial and
$V_{m}(x)\in W^{1, \infty}({\mathbb R})$. Hence, $V_{m}(x)$ and its derivative
$\displaystyle{\frac{dV_{m}}{dx}}$ are bounded on ${\mathbb R}$. Then
it can be trivially checked that the multiplication operators
\begin{equation}
\label{T}  
T_{m}\phi(x):=V_{m}(x)\phi(x), \quad \phi(x)\in H^{1}({\mathbb R})
\end{equation}
satisfy Assumption 1.1. 

\noindent
By virtue of the Sobolev inequality in one dimension (see e.g. Sect 8.5 of
~\cite{LL97}), we have
\begin{equation}
\label{e}
\|\phi(x)\|_{L^{\infty}({\mathbb R})}\leq \frac{1}{\sqrt{2}}
\|\phi(x)\|_{H^{1}({\mathbb R})}. 
\end{equation}  
Let us recall the algebra property for our Sobolev space.

\noindent
For any $\phi_{1}(x), \ \phi_{2}(x)\in H^{1}({\mathbb R})$, we have
$\phi_{1}(x)\phi_{2}(x)\in H^{1}({\mathbb R})$ and 
\begin{equation}
\label{alg}
\|\phi_{1}(x)\phi_{2}(x)\|_{H^{1}({\mathbb R})}\leq c_{a}\|\phi_{1}(x)\|_{H^{1}({\mathbb R})}
\|\phi_{2}(x)\|_{H^{1}({\mathbb R})},
\end{equation}
where $c_{a}>0$ is a constant.

\noindent
By means of the Young's inequality (see e.g.
Section 4.2 of ~\cite{LL97}), we obtain the upper bound on the norm of the
convolution as
\begin{equation}
\label{y}
\|\phi_{1}*\phi_{2}\|_{L^{2}({\mathbb R})}\leq \|\phi_{1}\|_{L^{1}({\mathbb R})}
\|\phi_{2}\|_{L^{2}({\mathbb R})}.
\end{equation}
Evidently, inequality (\ref{y}) implies the estimate
\begin{equation}
\label{yd}
\Big\|\frac{d}{dx}\int_{-\infty}^{\infty}\phi_{1}(x-y)\phi_{2}(y)dy\Big\|
_{L^{2}({\mathbb R})}\leq
\Big\|\frac{d\phi_{1}}{dx}\Big\|_{L^{1}({\mathbb R})}\|\phi_{2}\|_{L^{2}({\mathbb R})}.
\end{equation} 
Let us look for solutions of the nonlinear system of equations
(\ref{p}) of the form
\begin{equation}
\label{r}
u(x)=u_{0}(x)+u_{p}(x),
\end{equation}
with $u_{0}(x)$ given by (\ref{vf0}) and
\begin{equation}
\label{up}
u_{p}(x):=(u_{p, 1}(x), u_{p, 2}(x), ..., u_{p, N}(x))^{T}.  
\end{equation}  
Obviously, we derive the perturbative system
$$
u_{p, m}(x)=
$$
\begin{equation}
\label{pert}
[T_{m}(u_{0, m}(x)+u_{p, m}(x))]\int_{-\infty}^{\infty}
K_{m}(x-y)g_{m}(u_{0}(y)+u_{p}(y))dy, 
\end{equation}
with $1\leq m\leq N$.

\noindent
Let us use a closed ball in the Sobolev space
\begin{equation}
\label{b}
B_{\rho}:=\{u(x)\in H^{1}({\mathbb R}, {\mathbb R}^{N}) \ | \
\|u\|_{H^{1}({\mathbb R}, {\mathbb R}^{N})}\leq \rho \}, \quad 0<\rho\leq 1.
\end{equation}
We seek solutions of the system of equations (\ref{pert}) as fixed points of
the auxiliary nonlinear problem
\begin{equation}
\label{aux}
u_{m}(x)=[T_{m}(u_{0, m}(x)+v_{m}(x))] \int_{-\infty}^{\infty}
K_{m}(x-y)g_{m}(u_{0}(y)+v(y))dy,
\end{equation}
where $1\leq m\leq N$ in the ball (\ref{b}).

\noindent
Let us define the closed ball in the space of $N$ dimensions as
\begin{equation}
\label{i}
I:=\Big\{z\in {\mathbb R}^{N} \ | \ |z|_{{\mathbb R}^{N}}\leq \frac{1}{\sqrt{2}}
(\|u_{0}\|_{H^{1}({\mathbb R}, {\mathbb R}^{N})}+1)\Big\}
\end{equation}
along with the closed ball in the space of $C^{1}(I, {\mathbb R}^{N})$
vector-functions, namely
$$
D_{M}:=
$$
\begin{equation}
\label{M}
\{g(z):=(g_{1}(z), g_{2}(z), ..., g_{N}(z))\in C^{1}(I, {\mathbb R}^{N})
\ | \ \|g\|_{C^{1}(I, {\mathbb R}^{N})}\leq M \}, 
\end{equation}
with $M>0$.
In such context the norm
\begin{equation}
\label{gn1}
\|g\|_{C^{1}(I, {\mathbb R}^{N})}:=\sum_{m=1}^{N}\|g_{m}\|_{C^{1}(I)},
\end{equation}
\begin{equation}
\label{gn2}
\|g_{m}\|_{C^{1}(I)}:=\|g_{m}\|_{C(I)}+\sum_{n=1}^{N}\Big\|\frac{\partial g_{m}}
{\partial z_{n}}\Big\|_{C(I)},
\end{equation}  
where $\|\phi\|_{C(I)}:=\hbox{max}_{z\in I}|\phi(z)|$ for $\phi\in C(I)$.

\bigskip

\noindent
{\bf Assumption 1.2.} {\it Let $1\leq m\leq N$. Assume that
$g_{m}(z): {\mathbb R}^{N}\to {\mathbb R}$, so that
$g_{m}(0)=0$. It is also assumed that $g(z)|_{I}\in D_{M}$ and
it does not vanish identically in the ball $I$}.

\bigskip

\noindent
Let $\tau_{g}$ be the operator defined by the right side of system (\ref{aux}),
so that $u = \tau_{g}v$.
Our first main statement is as follows.

\bigskip

\noindent
{\bf Theorem 1.3.} {\it Let Assumptions 1.1 and 1.2 hold and
\begin{equation}
\label{rub}
c_{a}M(\|u_{0}\|_{H^{1}({\mathbb R}, {\mathbb R}^{N})}+1)^{2}Q\leq \frac{\rho}{2}.
\end{equation}
Then the map $\tau_{g}: B_{\rho}\to B_{\rho}$ associated with the system of
equations (\ref{aux}) is a strict contraction.
The unique fixed point $u_{p}(x)$ of this map $\tau_{g}$ is the only solution of
problem (\ref{pert}) in $B_{\rho}$.}

\bigskip

\noindent
Clearly, the resulting solution of system (\ref{p}) given by (\ref{r}) will
not vanish identically on the real line since all $g_{m}(0)=0$, the operators
$T_{m}$ are linear and the functions $u_{0, m}(x)$ are nontrivial for some
$1\leq m\leq N$ as assumed.
  
\bigskip

\noindent
For the technical purpose we introduce
\begin{equation}
\label{sig}
\sigma:=2c_{a}QM(\|u_{0}\|_{H^{1}({\mathbb R}, {\mathbb R}^{N})}+1)>0. 
\end{equation}
Our second major proposition is about the continuity of the cumulative solution
of problem (\ref{p}) given by formula (\ref{r}) with respect to the nonlinear
vector-function $g$.  

\bigskip

\noindent
{\bf Theorem 1.4.} {\it Suppose that the assumptions of Theorem 1.3 are valid.
Let $u_{p,j}(x)$ be the unique fixed point of the map
$\tau_{g_{j}}: B_{\rho}\to B_{\rho}$ for $g_{j}, \ j=1,2$ and the resulting solution
of problem (\ref{p}) with $g(z)=g_{j}(z)$ be given by
\begin{equation}
\label{u12}  
u_{j}(x)=u_{0}(x)+u_{p, j}(x).
\end{equation}
Then} 
\begin{equation}
\label{cont}
\|u_{1}(x)-u_{2}(x)\|_{H^{1}({\mathbb R}, {\mathbb R}^{N})}\leq
\end{equation}
$$
\frac{\sigma}
{2M (1-\sigma)}(\|u_{0}\|_{H^{1}({\mathbb R}, {\mathbb R}^{N})}+1)\|g_{1}(z)-g_{2}(z)\|_
{C^{1}(I, {\mathbb R}^{N})}.
$$

\bigskip

\noindent
We proceed to the proof of our first main result.

\bigskip 


\setcounter{equation}{0}

\section{The existence of the perturbed solution}

\noindent
{\it Proof of Theorem 1.3.} We choose an arbitrary $v(x)\in B_{\rho}$.
By virtue of (\ref{aux}) along with (\ref{alg}), we derive the estimate from
above for $1\leq m\leq N$ as
\begin{equation}
\label{auxub} 
\|u_{m}\|_{H^{1}({\mathbb R})}\leq
\end{equation}
$$
c_{a}\|T_{m}(u_{0, m}(x)+v_{m}(x))\|_{H^{1}({\mathbb R})}
\Big\|\int_{-\infty}^{\infty}K_{m}(x-y)g_{m}(u_{0}(y)+v(y))dy\Big\|_{H^{1}({\mathbb R})}.
$$
Let us consider the right side of (\ref{auxub}). Evidently,
\begin{equation}
\label{tov}  
\|T_{m}(u_{0, m}(x)+v_{m}(x))\|_{H^{1}({\mathbb R})}\leq \|T_{m}\|
(\|u_{0}(x)\|_{H^{1}({\mathbb R}, {\mathbb R}^{N})}+1). 
\end{equation}
By means of bound (\ref{y}), we have
$$
\Big\|\int_{-\infty}^{\infty}K_{m}(x-y)g_{m}(u_{0}(y)+v(y))dy\Big\|_{L^{2}({\mathbb R})}
\leq
$$
\begin{equation}
\label{kgvl2ub}
\|K_{m}\|_{L^{1}({\mathbb R})}\|g_{m}(u_{0}(x)+v(x))\|_{L^{2}({\mathbb R})}.
\end{equation}
Clearly, inequality (\ref{yd}) gives us
$$
\Big\|\frac{d}{dx}\int_{-\infty}^{\infty}K_{m}(x-y)g_{m}(u_{0}(y)+v(y))dy\Big\|
_{L^{2}({\mathbb R})}\leq
$$
\begin{equation}
\label{kgvl2ubd}
\Big\|\frac{dK_{m}}{dx}\Big\|_{L^{1}({\mathbb R})}\|g_{m}(u_{0}(x)+v(x))\|_
{L^{2}({\mathbb R})}.
\end{equation}
Formulas (\ref{kgvl2ub}) and (\ref{kgvl2ubd}) easily imply that
$$
\Big\|\int_{-\infty}^{\infty}K_{m}(x-y)g_{m}(u_{0}(y)+v(y))dy\Big\|_{H^{1}({\mathbb R})}
\leq
$$
\begin{equation}
\label{kgvl2ubc}
\|K_{m}\|_{W^{1, 1}({\mathbb R})}\|g_{m}(u_{0}(x)+v(x))\|_{L^{2}({\mathbb R})}.
\end{equation}
Note that, for $1\leq m\leq N$,
\begin{equation}
\label{gi}
g_{m}(u_{0}(x)+v(x))=\int_{0}^{1}
\nabla g_{m}(t(u_{0}(x)+v(x))).(u_{0}(x)+v(x))dt,
\end{equation}
where the dot stands for the scalar product of two vectors in ${\mathbb R}^{N}$.

\noindent
Let us recall inequality (\ref{e}). Hence, for $v(x)\in B_{\rho}$ we have
\begin{equation}
\label{u0pv}
|u_{0}+v|_{{\mathbb R}^{N}}\leq \frac{1}{\sqrt{2}}
(\|u_{0}\|_{H^{1}({\mathbb R}, {\mathbb R}^{N})}+1). 
\end{equation}
Thus,
$$
|g_{m}(u_{0}(x)+v(x))|\leq
$$
\begin{equation}
\label{G}
\hbox{sup}_{z\in I}|\nabla g_{m}(z)|_{{\mathbb R}^{N}}
|u_{0}(x)+v(x)|_{{\mathbb R}^{N}}\leq
M|u_{0}(x)+v(x)|_{{\mathbb R}^{N}},
\end{equation}
with the ball $I$ introduced in (\ref{i}), which yields
\begin{equation}
\label{Gn}
\|g_{m}(u_{0}(x)+v(x))\|_{L^{2}({\mathbb R})}
\leq M (\|u_{0}\|_{H^{1}({\mathbb R}, {\mathbb R}^{N})}+1).  
\end{equation}
Let us combine the estimates above to derive, for $1\leq m\leq N$,
\begin{equation}
\label{ur}  
\|u_{m}(x)\|_{H^{1}({\mathbb R})}\leq c_{a}M
(\|u_{0}\|_{H^{1}({\mathbb R}, {\mathbb R}^{N})}+1)^{2}
\|T_{m}\|\|K_{m}\|_{W^{1, 1}({\mathbb R})},
\end{equation}
such that
\begin{equation}
\label{urvf}  
\|u(x)\|_{H^{1}({\mathbb R}, {\mathbb R}^{N})}\leq c_{a}M
(\|u_{0}\|_{H^{1}({\mathbb R}, {\mathbb R}^{N})}+1)^{2}Q.
\end{equation}
By means of (\ref{rub}), we obtain that
$\|u(x)\|_{H^{1}({\mathbb R}, {\mathbb R}^{N})}\leq \rho$.
Hence, the vector-function $u(x)$, which is uniquely determined by (\ref{aux}) 
is contained in $B_{\rho}$ as well.
This means that the system of equations (\ref{aux}) defines a map
$\tau_{g}: B_{\rho}\to B_{\rho}$ under the stated assumptions.

\noindent
Our goal is to demonstrate that under the given conditions such map is a strict
contraction. Let us choose
arbitrary $v_{1,2}(x)\in B_{\rho}$. By virtue of the argument above,
$u_{1,2}:=\tau_{g}v_{1,2}\in B_{\rho}$. According to (\ref{aux}), for
$1\leq m\leq N$
$$
u_{1, m}(x)=
$$
\begin{equation}
\label{aux1}
[T_{m}(u_{0, m}(x)+v_{1, m}(x))]\int_{-\infty}^{\infty}
K_{m}(x-y)g_{m}(u_{0}(y)+v_{1}(y))dy,
\end{equation}
$$
u_{2, m}(x)=
$$
\begin{equation}
\label{aux2}
[T_{m}(u_{0, m}(x)+v_{2, m}(x))]\int_{-\infty}^{\infty}
K_{m}(x-y)g_{m}(u_{0}(y)+v_{2}(y))dy.
\end{equation}
By means of (\ref{aux1}) and (\ref{aux2}),
\begin{equation}
\label{u1mu2} 
u_{1, m}(x)-u_{2, m}(x)=
\end{equation}
$$
[T_{m}v_{1, m}(x)-T_{m}v_{2, m}(x)]\int_{-\infty}^{\infty}K_{m}(x-y)
g_{m}(u_{0}(y)+v_{1}(y))dy+
$$
$$
[T_{m}(u_{0, m}(x)+v_{2, m}(x))]\times
$$
$$
\int_{-\infty}^{\infty}K_{m}(x-y)
[g_{m}(u_{0}(y)+v_{1}(y))-g_{m}(u_{0}(y)+v_{2}(y))]dy.
$$
From (\ref{u1mu2}) using (\ref{alg}) we easily deduce that
\begin{equation}
\label{u1mu2hin}
\|u_{1, m}(x)-u_{2, m}(x)\|_{H^{1}({\mathbb R})}\leq c_{a}\|T_{m}v_{1, m}(x)-
T_{m}v_{2, m}(x)\|_{H^{1}({\mathbb R})}\times
\end{equation}
$$  
\Big\|\int_{-\infty}^{\infty}K_{m}(x-y)g_{m}(u_{0}(y)+v_{1}(y))dy\Big\|_
{H^{1}({\mathbb R})}+
$$
$$
c_{a}\|T_{m}(u_{0, m}(x)+v_{2, m}(x))\|_{H^{1}({\mathbb R})}\times
$$
$$
\Big\|\int_{-\infty}^{\infty}K_{m}(x-y)
[g_{m}(u_{0}(y)+v_{1}(y))-g_{m}(u_{0}(y)+v_{2}(y))]dy\Big\|_{H^{1}({\mathbb R})}.
$$
We derive the estimate from above on the right side of (\ref{u1mu2hin}). 
Clearly,
\begin{equation}
\label{tv1mtv2}  
\|T_{m}v_{1, m}(x)-T_{m}v_{2, m}(x)\|_{H^{1}({\mathbb R})}\leq \|T_{m}\|
\|v_{1}(x)-v_{2}(x)\|_{H^{1}({\mathbb R}, {\mathbb R}^{N})}.
\end{equation}
Note that by the same reasoning as above     
$$
\Big\|\int_{-\infty}^{\infty}K_{m}(x-y)g_{m}(u_{0}(y)+v_{1}(y))dy\Big\|_
{H^{1}({\mathbb R})}\leq
$$
\begin{equation}
\label{kgv1intl2h1}  
\|K_{m}\|_{W^{1, 1}({\mathbb R})}M(\|u_{0}\|_{H^{1}({\mathbb R}, {\mathbb R}^{N})}+1).
\end{equation}
Hence, the first term in the right side of formula (\ref{u1mu2hin}) can
be bounded from above by
\begin{equation}
\label{catv12km}
c_{a}\|T_{m}\|\|v_{1}(x)-v_{2}(x)\|_{H^{1}({\mathbb R}, {\mathbb R}^{N})}
\|K_{m}\|_{W^{1, 1}({\mathbb R})}M
(\|u_{0}\|_{H^{1}({\mathbb R}, {\mathbb R}^{N})}+1).
\end{equation}
Thus, it remains to consider the second term in the right side of
(\ref{u1mu2hin}). Obviously, for $1\leq m\leq N$,
\begin{equation}
\label{tu0v2h1}
\|T_{m}(u_{0, m}(x)+v_{2, m}(x))\|_{H^{1}({\mathbb R})}\leq \|T_{m}\|
(\|u_{0}\|_{H^{1}({\mathbb R}, {\mathbb R}^{N})}+1).  
\end{equation}
Using (\ref{y}) and (\ref{yd}), we easily derive
$$
\Big\|\int_{-\infty}^{\infty}K_{m}(x-y)[g_{m}(u_{0}(y)+v_{1}(y))-
g_{m}(u_{0}(y)+v_{2}(y))]dy\Big\|_{H^{1}({\mathbb R})}\leq
$$
\begin{equation}
\label{kgu0v12l2h1}
\|K_{m}\|_{W^{1, 1}({\mathbb R})}\|
g_{m}(u_{0}(x)+v_{1}(x))-g_{m}(u_{0}(x)+v_{2}(x))\|_{L^{2}({\mathbb R})}.      
\end{equation}
Let us express, for $1\leq m\leq N$,
$$
g_{m}(u_{0}(x)+v_{1}(x))-g_{m}(u_{0}(x)+v_{2}(x))=
$$
\begin{equation}
\label{gu0v12iz}
\int_{0}^{1}\nabla
g_{m}(u_{0}(x)+tv_{1}(x)+(1-t)v_{2}(x)).[v_{1}(x)-v_{2}(x)]dt.
\end{equation}  
Formula (\ref{gu0v12iz}) implies that
\begin{equation}
\label{gu0v12m}
|g_{m}(u_{0}(x)+v_{1}(x))-g_{m}(u_{0}(x)+v_{2}(x))|\leq
M|v_{1}(x)-v_{2}(x)|_{{\mathbb R}^{N}},
\end{equation}
so that for the norm
$$
\|g_{m}(u_{0}(x)+v_{1}(x))-g_{m}(u_{0}(x)+v_{2}(x))\|_{L^{2}({\mathbb R})}\leq
$$
\begin{equation}
\label{gu0v12n}
M\|v_{1}(x)-v_{2}(x)\|_{H^{1}({\mathbb R}, {\mathbb R}^{N})}.
\end{equation}
Hence, the second term in the right side of bound (\ref{u1mu2hin}) can
be estimated from above by expression (\ref{catv12km}) as well.

\noindent
We arrive at
$$
\|u_{1}(x)-u_{2}(x)\|_{H^{1}({\mathbb R}, {\mathbb R}^{N})}\leq
$$
\begin{equation}
\label{u12v12h1}  
2c_{a}QM(\|u_{0}\|_{H^{1}({\mathbb R}, {\mathbb R}^{N})}+1)
\|v_{1}(x)-v_{2}(x)\|_{H^{1}({\mathbb R}, {\mathbb R}^{N})}.
\end{equation}
Formula (\ref{u12v12h1}) and definition (\ref{sig}) give us
\begin{equation}
\label{tgv12}
\|\tau_{g}v_{1}(x)-\tau_{g}v_{2}(x)\|_{H^{1}({\mathbb R}, {\mathbb R}^{N})}\leq \sigma
\|v_{1}(x)-v_{2}(x)\|_{H^{1}({\mathbb R}, {\mathbb R}^{N})}.
\end{equation}
It can be trivially checked via (\ref{rub}) that the constant in the right
side of inequality above
\begin{equation}
\label{sigm}  
\sigma<1.
\end{equation}
Therefore, the map $\tau_{g}: B_{\rho}\to B_{\rho}$ defined
by system (\ref{aux}) is a strict contraction under the stated assumptions.
Its unique fixed point $u_{p}(x)$ is the only solution of the system of
equations (\ref{pert}) in the ball $B_{\rho}$. 
The resulting $u(x)$ given by (\ref{r}) solves problem (\ref{p}).
\hspace{6cm} $\Box$

\bigskip

\noindent
We conclude the article by demonstrating the validity of the second main
statement.

\bigskip


\setcounter{equation}{0}

\section{The continuity of the resulting solution with respect to
the vector function $g$}

\noindent
{\it Proof of Theorem 1.4.} Evidently, under the given conditions
\begin{equation}
\label{up1up2}  
u_{p,1}=\tau_{g_{1}}u_{p,1}, \quad u_{p,2}=\tau_{g_{2}}u_{p,2}.
\end{equation}
Hence,
\begin{equation}
\label{up1mup2}  
u_{p,1}-u_{p,2}=\tau_{g_{1}}u_{p,1}-\tau_{g_{1}}u_{p,2}+\tau_{g_{1}}u_{p,2}-
\tau_{g_{2}}u_{p,2},
\end{equation}
such that
$$
\|u_{p,1}-u_{p,2}\|_{H^{1}({\mathbb R}, {\mathbb R}^{N})}\leq
$$
\begin{equation}
\label{up1mup2n}  
\|\tau_{g_{1}}u_{p,1}-\tau_{g_{1}}u_{p,2}\|_{H^{1}({\mathbb R}, {\mathbb R}^{N})}+
\|\tau_{g_{1}}u_{p,2}-\tau_{g_{2}}u_{p,2}\|_{H^{1}({\mathbb R}, {\mathbb R}^{N})}.
\end{equation}
By virtue of inequality (\ref{tgv12}), we have
\begin{equation}
\label{tg1up1up2}  
\|\tau_{g_{1}}u_{p,1}-\tau_{g_{1}}u_{p,2}\|_{H^{1}({\mathbb R}, {\mathbb R}^{N})}\leq 
\sigma\|u_{p,1}-u_{p,2}\|_{H^{1}({\mathbb R}, {\mathbb R}^{N})},
\end{equation}
where $\sigma$ is defined in (\ref{sig}), so that (\ref{sigm}) is valid.

\noindent
Thus, we arrive at
\begin{equation}
\label{sigma}
(1-\sigma)\|u_{p,1}-u_{p,2}\|_{H^{1}({\mathbb R}, {\mathbb R}^{N})}\leq
\|\tau_{g_{1}}u_{p,2}-\tau_{g_{2}}u_{p,2}\|_{H^{1}({\mathbb R}, {\mathbb R}^{N})}.
\end{equation}
Clearly, for our fixed point $\tau_{g_{2}}u_{p,2}=u_{p,2}$. Let us introduce
$\eta(x):=\tau_{g_{1}}u_{p,2}$, such that for $1\leq m\leq N$
$$
\eta_{m}(x)=
$$
\begin{equation}
\label{12}
[T_{m}(u_{0,m}(x)+u_{p,2,m}(x))]\int_{-\infty}^{\infty}
K_{m}(x-y)g_{1,m}(u_{0}(y)+u_{p,2}(y))dy,
\end{equation}
$$
u_{p,2,m}(x)=
$$
\begin{equation}
\label{22}
[T_{m}(u_{0,m}(x)+u_{p,2,m}(x))]\int_{-\infty}^{\infty}
K_{m}(x-y)g_{2,m}(u_{0}(y)+u_{p,2}(y))dy.
\end{equation}
By means of formulas (\ref{12}) and (\ref{22}) 
$$
\eta_{m}(x)-u_{p,2,m}(x)=[T_{m}(u_{0,m}(x)+u_{p,2,m}(x))]\times
$$
\begin{equation}
\label{rup2}  
\int_{-\infty}^{\infty}K_{m}(x-y)
[g_{1,m}(u_{0}(y)+u_{p,2}(y))-g_{2,m}(u_{0}(y)+u_{p,2}(y))]dy.
\end{equation}
Let us recall bound (\ref{alg}). Hence,
$$
\|\eta_{m}(x)-u_{p,2,m}(x)\|_{H^{1}({\mathbb R})}\leq c_{a}
\|T_{m}(u_{0,m}(x)+u_{p,2,m}(x))\|_{H^{1}({\mathbb R})}\times
$$
\begin{equation}
\label{rup2n}
\Big\|\int_{-\infty}^{\infty}K_{m}(x-y)
[g_{1,m}(u_{0}(y)+u_{p,2}(y))-g_{2,m}(u_{0}(y)+u_{p,2}(y))]dy\Big\|_{H^{1}({\mathbb R})}.
\end{equation}
Obviously, the estimate from above
\begin{equation}
\label{tu0up2n}  
\|T_{m}(u_{0,m}(x)+u_{p,2,m}(x))\|_{H^{1}({\mathbb R})}\leq \|T_{m}\|
(\|u_{0}\|_{H^{1}({\mathbb R}, {\mathbb R}^{N})}+1)
\end{equation}
holds.

\noindent
By virtue of (\ref{y}) and (\ref{yd}), we easily obtain
$$
\Big\|\int_{-\infty}^{\infty}K_{m}(x-y)
[g_{1,m}(u_{0}(y)+u_{p,2}(y))-g_{2,m}(u_{0}(y)+u_{p,2}(y))]dy\Big\|_{H^{1}({\mathbb R})}
\leq
$$
\begin{equation}
\label{kg1g2u0up2nh}  
\|K_{m}\|_{W^{1, 1}({\mathbb R})}\|g_{1,m}(u_{0}(x)+u_{p,2}(x))-
g_{2,m}(u_{0}(x)+u_{p,2}(x))\|_{L^{2}({\mathbb R})}.
\end{equation}
Evidently, for $1\leq m\leq N$,
\begin{equation}
\label{g12u0up2}  
[g_{1,m}-g_{2,m}](u_{0}(x)+u_{p,2}(x))=
\end{equation}
$$
\int_{0}^{1}\nabla[g_{1,m}-g_{2,m}](t(u_{0}(x)+u_{p,2}(x))).[u_{0}(x)+u_{p,2}(x)]dt.
$$
Let us use (\ref{g12u0up2}) to derive
\begin{equation}
\label{g12u0up2m}  
|[g_{1,m}-g_{2,m}](u_{0}(x)+u_{p,2}(x))|\leq
\|g_{1,m}-g_{2,m}\|_{C^{1}(I)}|u_{0}(x)+u_{p,2}(x)|_{{\mathbb R}^{N}}.
\end{equation}
Then for the norm
$$
\|[g_{1,m}-g_{2,m}](u_{0}(x)+u_{p,2}(x))\|_{L^{2}({\mathbb R})}\leq
$$
\begin{equation}
\label{g12u0up2mn}
\|g_{1}-g_{2}\|_{C^{1}(I, {\mathbb R}^{N})}(\|u_{0}\|_{H^{1}({\mathbb R}, {\mathbb R}^{N})}+1).
\end{equation}
By means of estimates (\ref{rup2n}), (\ref{tu0up2n}), 
(\ref{kg1g2u0up2nh}), (\ref{g12u0up2mn}),  we arrive at
$$
\|\eta(x)-u_{p,2}(x)\|_{H^{1}({\mathbb R}, {\mathbb R}^{N})}\leq
$$
\begin{equation}
\label{rup2h1r}
c_{a}Q
(\|u_{0}\|_{H^{1}({\mathbb R}, {\mathbb R}^{N})}+1)^{2}\|g_{1}(z)-g_{2}(z)\|_
{C^{1}(I, {\mathbb R}^{N})}.
\end{equation}
Formulas (\ref{sigma}) and (\ref{rup2h1r}) imply that
$$
\|u_{p,1}(x)-u_{p,2}(x)\|_{H^{1}({\mathbb R}, {\mathbb R}^{N})}\leq
$$
\begin{equation}
\label{p1p2}
\frac{c_{a}}{1-\sigma}
(\|u_{0}\|_{H^{1}({\mathbb R}, {\mathbb R}^{N})}+1)^{2}Q
\|g_{1}(z)-g_{2}(z)\|_{C^{1}(I, {\mathbb R}^{N})}.
\end{equation}
By virtue of (\ref{u12}) along with inequality (\ref{p1p2}) and definition
(\ref{sig}), bound (\ref{cont}) is valid.  \hspace{9.5cm} $\Box$

\bigskip

\noindent
{\bf Remark 3.1.} {\it The results of the present work will be generalized to
higher dimensions in forthcoming articles.}      

\bigskip


\section{Acknowledgement} The work was partially supported by the
NSERC Discovery grant.

\bigskip


\end{document}